\documentclass[Afour,sageh,times]{sagej}

\pdfoutput=1

\usepackage{moreverb,url}

\usepackage[pdftex,rgb]{xcolor}
\usepackage[colorlinks,bookmarksopen,bookmarksnumbered,citecolor=red,urlcolor=red]{hyperref}

\newcommand\BibTeX{{\rmfamily B\kern-.05em \textsc{i\kern-.025em b}\kern-.08em
T\kern-.1667em\lower.7ex\hbox{E}\kern-.125emX}}

\usepackage{algpseudocode}
\usepackage{algorithm}

\usepackage{listings}
\usepackage{colortbl}
\usepackage{multirow, bigstrut}
\usepackage{booktabs}
\usepackage{siunitx}

\definecolor{tab-blue}{RGB}{31,119,180}
\definecolor{tab-orange}{RGB}{255,127,14}
\definecolor{tab-green}{RGB}{44,160,44}
\definecolor{tab-red}{RGB}{214,39,40}
\definecolor{tab-purple}{RGB}{148,103,189}
\definecolor{tab-brown}{RGB}{140,86,75}
\definecolor{tab-pink}{RGB}{227,119,194}
\definecolor{tab-gray}{RGB}{127,127,127}
\definecolor{tab-olive}{RGB}{188,189,34}
\definecolor{tab-cyan}{RGB}{23,190,207}

\usepackage[export]{adjustbox}

\renewcommand{\vec}{\mathbf}
\newcommand{\mytext}[2]{{\fontfamily{lmss}\selectfont\textbf{\textcolor{#1}{#2}}}}

\usepackage{subcaption}
\usepackage[capitalise,noabbrev]{cleveref}

\begin{document}

\runninghead{Spies, Olson, MacLachlan}

\title{Exploiting mesh structure to improve multigrid performance for saddle point problems}

\author{Lukas Spies\affilnum{1}, Luke Olson\affilnum{1}, and Scott MacLachlan\affilnum{2}}

\affiliation{\affilnum{1}University of Illinois Urbana-Champaign, USA\\
\affilnum{2}Memorial University of Newfoundland, Canada}

\corrauth{Lukas Spies,
Siebel Center for Computer Science,
University of Illinois Urbana-Champaign,
201 N. Goodwin Ave,
Urbana, IL 61801,
USA}

\email{lspies@illinois.edu}

\begin{abstract}
     In recent years, solvers for finite-element discretizations of linear or linearized saddle-point problems, like the Stokes and Oseen equations, have become well established. There are two main classes of preconditioners for such systems: those based on block-factorization approach and those based on monolithic multigrid. Both classes of preconditioners have several critical choices to be made in their composition, such as the selection of a suitable relaxation scheme for monolithic multigrid. From existing studies, some insight can be gained as to what options are preferable in low-performance computing settings, but there are very few fair comparisons of these approaches in the literature, particularly for modern architectures, such as GPUs.
     In this paper, we perform a comparison between a block-triangular preconditioner and a monolithic multigrid method with the three most common choices of relaxation scheme - Braess-Sarazin, Vanka, and Schur-Uzawa. We develop a performant Vanka relaxation algorithm for structured-grid discretizations, which takes advantage of memory efficiencies in this setting.
     We detail the behavior of the various CUDA kernels for the multigrid relaxation schemes and evaluate their individual arithmetic intensity, performance, and runtime. Running a preconditioned FGMRES solver for the Stokes equations with these preconditioners allows us to compare their efficiency in a practical setting. We show that
monolithic multigrid can outperform block-triangular preconditioning, and that using Vanka or Braess-Sarazin relaxation is most efficient.  Even though multigrid with Vanka relaxation exhibits reduced performance on the CPU (up to $100\%$ slower than Braess-Sarazin), it is able to outperform Braess-Sarazin by more than $20\%$ on the GPU, making it a competitive algorithm, especially given the high amount of algorithmic tuning needed for effective Braess-Sarazin relaxation.
\end{abstract}


\maketitle

\section{Introduction}

    Finite-element discretizations are a popular choice for coupled systems such as magnetohydrodynamics (MHD), or the Stokes or Navier-Stokes equations. Even though solvers for finite-element discretizations of such saddle-point problems are well established, designing efficient and scalable solvers on emerging computing architectures for such systems remains an ongoing challenge~\citep{SolverCurrent1,SolverCurrent2,SolverCurrent3}.

    Here, we focus on preconditioned Krylov methods for the linear or linearized systems of equations that arise in solving such problems.  There are two main classes of preconditioners for such systems: preconditioners based on block-factorization approaches~\cite{HCElman_DJSilvester_AWathen_2005a,MBenzi_GHGolub_JLiesen_2005a,SolverBlock} and those based on monolithic multigrid principles~\citep{MultigridStokes, LowOrderPreconditioning}. Within each class there is considerable variability in their building blocks, such as the choice of relaxation scheme in monolithic multigrid, including Braess-Sarazin~\citep{BraessSarazin, WZulehner_2000a}, Vanka~\citep{Vanka2}, and Schur-Uzawa~\citep{JFMaitre_FMusy_PNignon_1984a} relaxation.

    From existing studies~\citep{BensonAdlerCyrMacLachlanTuminaro,FourierVankaStokes,LowOrderPreconditioning,JAdler_etal_2023a} some insight can be gained into which algorithms are preferable in serial (low-performance) computing settings, but there are relatively few fair comparisons of these approaches in literature~\citep{VKBS1,VKBS2,CompVankaBSStokes,JAdler_etal_2015b}, in particular for geometric multigrid on modern architectures, such as GPUs.  This is particularly important given changes in prevailing HPC architectures in the past two decades. Of interest to us is that Vanka relaxation has been shown to lead to scalable performance mathematically~\citep{VankaSmoothing,FourierVankaStokes} and is seen as an algorithm that is well-suited for implementation on modern GPUs but, to our knowledge, no performance studies support this claim. As we discuss below, one main difficulty in getting good performance out of the Vanka algorithm is the high cost of memory movement for forming the various Vanka patches and updating the global solution, which requires a careful approach to achieve good performance. Similar issues have recently been considered using related additive Schwarz relaxation schemes within multigrid applied to the Poisson equation~\citep{MunchKronbichlerAS}, where it was found that, with optimization of memory caching and reducing communication between patches, additive Schwarz method built around cell-centric patches are capable of outperforming optimized point-Jacobi-based relaxation schemes.

    More broadly, the parallel scalability of multigrid algorithms on modern architectures faces many challenges related to indirection and increased coarse-grid communication costs~\citep{2016_BiFaGrOlSc,2020_BiGrOl_reducing}. This makes data locality and the cost of data movement a central issue, but one that can be solved by carefully exploiting structure in the problem as is done, e.g., in black box multigrid (BoxMG) algorithms~\citep{JEDendy_1982a,ScalingStructuredMG,ScalableLineRelaxation}. In this work, we use highly structured meshes that allow us to encode various information about the data and how it is accessed in the structure itself, similarly to the BoxMG paradigm. Notably, we work with a structured matrix representation and implement algorithms that take full advantage of this, in order to minimize memory accesses and maximize floating point operations (arithmetic intensity).

    In this paper, we first introduce the Stokes equations as our model problem and provide an overview of their structure and resulting discretization. We then introduce two different preconditioners for the FGMRES algorithm used to solve such equations, monolithic multigrid and the upper block-triangular preconditioner. For monolithic multigrid, we introduce three different relaxation schemes, Braess-Sarazin, Vanka, and Schur-Uzawa. As Braess-Sarazin and Vanka are two common choices for relaxation schemes, we focus our performance analysis on these, noting that Schur-Uzawa can be implemented with the same kernels as Braess-Sarazin. After illustrating which kernels are the biggest contributors to each algorithm, we then break the performance analysis into two parts: common kernels (matrix-vector and array operations) and Vanka-specific kernels (forming patches, updating the global solution, solving patch systems). For each part, the arithmetic intensity, performance, and runtime are analyzed in order to gain a full understanding of the algorithms and how they compare. This leads to a roofline model to investigate how much better the kernels could be doing, if at all. To finish our performance analysis, we compare a full solve of the Stokes equations using FGMRES preconditioned with both a block-triangular preconditioner and a multigrid V-cycle preconditioner with Braess-Sarazin, Vanka, and Schur-Uzawa as relaxation schemes. We show that Vanka is, indeed, a competitive algorithm on modern architectures with careful design. We also show that simply porting a performant CPU implementation of Vanka directly to the GPU is not sufficient for achieving competitive performance.

\section{The Stokes equations and their discretization}

    \subsection{Problem setup}

        Fluid flow where viscous forces are much greater than advective inertia is called Stokes flow. In nature, flow with such properties occurs in many places, e.g., in geodynamics or in the swimming movement of microorganisms. The equations of motions arising from this flow are called the Stokes equations and can be viewed as a simplification of the steady-state Navier-Stokes equations in the limit of small Reynolds number, $Re \ll 1$. They are not only well-suited to be solved by iterative solvers~\citep{StokesIterative}, but they also serve as a suitable prototype for a wide range of models that lead to saddle-point structure.

\begin{figure}[ht]
            \centering
            \includegraphics{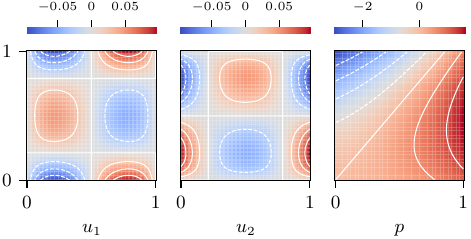}
            \caption{Visualization of three components of the manufactured solution.}\label{fig:samplesol}
        \end{figure}

        Specifically, we consider the incompressible Stokes equation with constant viscosity $\nu$ in the unit-square domain $\Omega = [0,1]^2 \in \mathbb{R}^2$. The equations are given by
        \begin{align}
            -\nu\nabla^2\vec{u} + \nabla p &= \vec{f},\label{eq:stokes1}\\
            \nabla\cdot \vec u &= 0.\label{eq:stokes2}
        \end{align}
        Dirichlet boundary conditions on velocity are enforced on all edges of the domain, but no boundary conditions are imposed on pressure. Here, we consider
        no-flux boundary conditions,
        \begin{equation}
            \vec u\cdot \vec n = 0 \;\text{on}\; \partial\Omega,
        \end{equation}
        where $\vec n$ is the outward pointing normal vector. Thus, we define the Hilbert space $\vec H_0^1(\Omega)$ as
        \begin{equation}
            \vec H_0^1(\Omega) = \{\vec v \in \vec H^1(\Omega) : \vec v \cdot \vec n = 0\;\text{on}\;\partial\Omega\}
        \end{equation}
        and $L^2(\Omega)/\mathbb{R}$ as the quotient space of equivalence classes of functions in $L^2(\Omega)$ differing by a constant. The weak form is then defined as: Find $(\vec u, p) \in H_0^1(\Omega) \times L^2(\Omega)/\mathbb{R}$ such that
        \begin{alignat}{2}
            a(\vec u, \vec v) + b(\vec v, p) &= (\vec f, \vec v) &&\quad\forall \vec v \in \vec H_0^1(\Omega)\label{eq:weak1}\\
            b(\vec u, q) &= 0 &&\quad\forall q \in L^2(\Omega)/\mathbb{R}\label{eq:weak2}
        \end{alignat}
        where
        \begin{align}
            a(\vec u, \vec v) &= \nu \int_\Omega \nabla \vec u : \nabla \vec v,\\
            b(\vec v, p) &= \int_\Omega p\nabla\cdot \vec v.
        \end{align}
        We work with a manufactured solution~\citep{BAyuso_etal_2014a, JAdler_etal_2015b} that satisfies the properties above, given by
        \begin{align}
            \vec u(x,y) &= \begin{cases} x(1-x)(2x-1)(6y^2-6y+1)\\ y(y-1)(2y-1)(6x^2-6x+1)\end{cases}\\
            p(x,y) &= x^2 - 3y^2 + \frac{8}{3}xy,
        \end{align}
        with $\vec f$ computed to satisfy~\eqref{eq:stokes1}.
        A visualization of this sample is found in~\cref{fig:samplesol}.

    \subsection{Discretization}

\begin{figure}[ht]
            \centering
            \begin{subfigure}[b]{0.48\columnwidth}
                \includegraphics{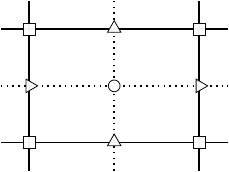}
                \caption{Q2 element}
            \end{subfigure}
            \begin{subfigure}[b]{0.48\columnwidth}
                \includegraphics{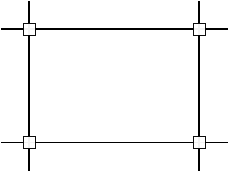}
                \caption{Q1 element}
            \end{subfigure}
            \caption{Illustration of the degrees of freedom for a Q2 and Q1 element, with different types of degrees of freedom identified by different shapes.}\label{fig:finele}
        \end{figure}

        We consider the standard Q2--Q1 Taylor-Hood mixed finite-element discretization on a uniform grid for discretizing the system in~\cref{eq:stokes1,eq:stokes2}. For the velocity, this uses Q2 elements, with biquadratic polynomials for each component on each element as a basis.  For the pressure, this uses Q1 elements, with bilinear polynomials on each element as a basis. Both velocities and pressures are required to be continuous across element boundaries. An illustration of the degrees of freedom in these elements is shown in~\cref{fig:finele}.
        The resulting Q2 and Q1 elements have a total of nine and four degrees of freedom per element, respectively. Such a discretization of the Stokes equations directly relates back to the weak form as shown in~\cref{eq:weak1,eq:weak2}, defining the matrices $L$ and $B$ by
        \begin{align}
        L_{i,j} & = a(\vec{\psi}_j,\vec{\psi}_i) \\
        B_{k,j} & = b(\vec{\psi}_j,\phi_k).
        \end{align}
        It is important to note that the introduction of a basis gives us two ``views'' on the finite-element approximations, writing $\vec{u} = \sum_i u_i\vec{\psi}_i$ and $p = \sum_k p_k\phi_k$, so we can consider the functions $\vec{u}$ and $p$ directly, or think about their coefficients in the basis expansion, $\{u_i\}$ and $\{p_k\}$.  In what follows, we follow the usual convention of overloading the notation $\vec{u}$ and $p$ to denote both the functions themselves and the vectors of basis coefficients, with the distinction typically clear from context, in that $L\vec{u}$ is the matrix acting on the basis coefficients, while $a(\vec{u},\vec{v})$ is the bilinear form evaluated on the function.  With this matrix representation, the solution of the weak form can be expressed as that of the linear system
        \begin{equation}
            \begin{bmatrix}
                L & B^T\\
                B & 0
            \end{bmatrix}
            \begin{bmatrix}
                \vec u\\
                p
            \end{bmatrix}
            =
            \begin{bmatrix}
                \vec f\\0
            \end{bmatrix}\label{eq:stokesmatrix}
        \end{equation}
        where $L$ is the discretized Laplacian, $B$ and $B^T$ are the discretized divergence of $\vec u$ and gradient of $p$, respectively, $\vec u$ and $p$ are the discretized velocity and pressure components of the solution, and $\vec f$ is the velocity component of the right-hand side.  In what follows, to save space, we will write the $2\times 2$ block matrix in~\cref{eq:stokesmatrix} as $A$.

        Such a system is challenging to solve, as it is symmetric but indefinite due to the zero block in the lower right-hand corner of the system matrix. This causes many common iterative methods (e.g., stationary methods like Jacobi and Gauss-Seidel) to not work as they typically involve inverting the diagonal of the system matrix.  Study of numerical methods for solution of such saddle-point systems is a well-established discipline~\citep{MBenzi_GHGolub_JLiesen_2005a}. One possible solution to these challenges is the use of block preconditioners, based on the block LU factorization of the coupled system.  This approach has been well-developed for discretizations of the Stokes equations~\citep{MBenzi_GHGolub_JLiesen_2005a,HCElman_DJSilvester_AWathen_2005a}.  However, existing studies~\citep{JAdler_etal_2015b} suggest that monolithic preconditioners can be more efficient. Thus, we also consider monolithic multigrid as preconditioner for FGMRES\@.  We note that monolithic multigrid preconditioners are generally not symmetric and positive definite and, so, they cannot be used directly as preconditioners for MINRES; however, the added computational work for orthogonalization in FGMRES is more than made up for by the quick convergence of the monolithic multigrid approach.

    \subsection{Structured matrix representation}

\begin{figure}[ht]
            \centering
            \includegraphics{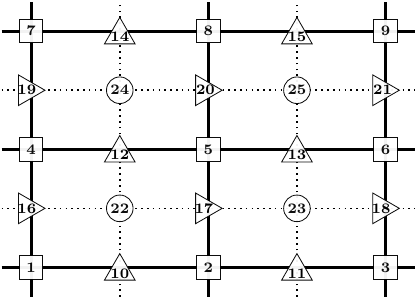}
            \caption{Local numbering of degrees of freedom around nodal degree of freedom $5$.}\label{fig:nodaldof}
        \end{figure}

        Any iterative solver naturally depends on calculations of matrix-vector products for the block-structured matrix in~\cref{eq:stokesmatrix}.  In general, such calculations require indirect addressing, when arbitrary numbers of elements can be adjacent to each node of the mesh, leading to irregular communication patterns.  However, when we restrict the mesh to have logically rectangular structure (meaning that each node is at the intersection of four edges, and is adjacent to four elements), then applying the discretization matrix can be done in a stencil-wise fashion, where each degree of freedom requires information from at most a $2\times2$ element patch. This allows storing the system matrices in a highly efficient data structure by numbering the degrees of freedom in lexicographic order.  For the Q1 discretization, this ordering is natural, since the only degrees of freedom occur at the nodes in the mesh, that can be labeled lexicographically by $(x,y)$-indices.  For the Q2 discretization, we separate the degrees of freedom into four sets, given by those associated with the nodes of the mesh, then those at midpoints of the $x$ and $y$ edges, and, finally, those associated the cell centers. \Cref{fig:nodaldof} shows a local numbering of those degrees of freedom on the $2\times 2$ patch around a node.

        With this ordering, the system matrix is stored as an array of arrays, where the first (outer) array index corresponds to the row in the matrix associated with that degree of freedom, and the entries inside that array correspond to the the degrees of freedom surrounding it, according to the numbering in~\cref{fig:nodaldof}. Exploiting the structure in this way allows us to minimize the memory storage needed and maximize the performance, as only a small amount of memory needs to be loaded and read for any operation. Additionally, no explicit pointers or indices need to be stored (as would be needed for standard sparse matrix representations), as this information is encoded in the data structure itself.

        This approach translates directly to higher dimensions. For a 3D discretization of this type on ``brick'' elements, we would extend the above to $2\times 2\times 2$ element patches with $5\times5\times5=125$ degrees of freedom ($5$ in each dimension).  These can be labeled in an analogous way and, consequently, stored similarly in a single array of arrays, with outer index corresponding to the rows in the system matrix, and inner index corresponding to the local numbering around each degree of freedom.

\section{Multigrid}

    Multigrid methods are based on the notion that standard (but slow-to-converge) iterative methods are generally effective at reducing some errors in a discrete approximation, but that the subspace of slow-to-converge errors of such an iteration is better treated by a complementary process~\citep{WLBriggs_VEHenson_SFMcCormick_2000a,UTrottenberg_etal_2001a}.  A natural approach to reduce the slow-to-converge errors is with correction from a coarser-grid realization of the same discretized problem, where those modes can be accurately resolved by recursively applying the same iterations to problems with fewer degrees of freedom, until some suitably coarse version of the problem is found where a sparse direct solver can be effectively applied.  For higher-order discretizations and systems of PDEs, such an error classification breaks down~\citep{YHe_SMacLachlan_2018a}, but the multigrid principle remains effective, in that we can define relaxation schemes that effectively damp a large portion of the error in a given approximation, and the remaining error can be effectively corrected from a coarse grid.

    The standard multigrid solution algorithm is known as the V-cycle, since it traverses a given hierarchy of meshes from the given finest grid to the coarsest, then back to the finest.  In the ``downward'' sweep of the traversal (from fine-to-coarse), on each level, an initial approximation (generally a zero vector) is improved by a specified relaxation scheme.  Then, the residual associated with that improved approximation is calculated and restricted to the next coarsest grid, where the scheme recurses.  On the ``upward'' sweep, the current approximation is improved by interpolating a correction back from the next coarser grid, then running relaxation again, before proceeding to the next finer grid.  For transferring residuals and corrections between grids, we define a single interpolation operator that maps from a coarse grid to the next finer grid, and use its transpose as a restriction operator.  In this work, we follow the standard geometric multigrid approach of using the finite-element interpolation operators, that naturally map from coarse-grid versions of the Q2 and Q1 spaces to their fine-grid analogues. \Cref{alg:mg} presents an algorithmic overview of the multigrid V-cycle for the Stokes equations, following the convention that level 0 is the coarsest grid in the hierarchy, and we are interested in the solution on some given fine grid, for fixed $l>0$.
    \begin{algorithm}[ht]
        \begin{algorithmic}[1]
            \Function{MG}{$A_l$, $\vec u_l$, $p_l$, $\vec f_l$, $g_l$, $l$}
                \State{Relax on $\vec u_l$ and $p_l$}
                \State{Compute residual: $\begin{bmatrix}\vec r_{\vec u,l}\\r_{p,l} \end{bmatrix} = \begin{bmatrix}\vec f_l\\ g_l \end{bmatrix} - A_l\begin{bmatrix} \vec u_l\\p_l \end{bmatrix}$}
                \State{Restriction: $\begin{bmatrix}\vec r_{\vec u,{l-1}}\\r_{p,{l-1}} \end{bmatrix} = P^T_{l-1}\begin{bmatrix}\vec r_{\vec u,l}\\r_{p,l} \end{bmatrix}$}
                \If{$l\text{ is }1$}
                    \State{$\begin{bmatrix}\vec e_{\vec u,0}\\e_{p,0} \end{bmatrix} = A_0^{-1}\begin{bmatrix}\vec r_{\vec u,0}\\r_{p,0} \end{bmatrix}$}
                \Else
                    \State{$\begin{bmatrix}\vec e_{\vec u,l-1}\\e_{p,l-1} \end{bmatrix} = \text{MG(}A_{l-1}, \vec 0, 0, \vec r_{\vec u,l-1},r_{p,l-1}, l-\nolinebreak1\text{)}$}
                \EndIf
                \State{Correction: $\begin{bmatrix} \vec u_l\\p_l \end{bmatrix} = \begin{bmatrix} \vec u_l\\p_l \end{bmatrix} + P_{l-1}\begin{bmatrix}\vec e_{\vec u,l-1}\\e_{p,l-1} \end{bmatrix}$}
                \State{Relax on $\vec u_l$ and $p_l$}
            \EndFunction
        \end{algorithmic}
        \caption{Multigrid V-cycle for Stokes equations}\label{alg:mg}
    \end{algorithm}

    While traditional relaxation schemes, such as (weighted) Jacobi or Gauss-Seidel are effective for elliptic problems, they generally cannot be applied directly to saddle-point systems, due to the zero block in the matrix.  Thus, specialized relaxation schemes are commonly developed and analyzed for the Stokes equations. In this work, we consider four different preconditioning approaches, comparing monolithic multigrid with Braess-Sarazin~\citep{BraessSarazin, WZulehner_2000a}, Vanka~\citep{Vanka2}, and Schur-Uzawa~\citep{JFMaitre_FMusy_PNignon_1984a} relaxation with an upper Block-Triangular preconditioner. We focus in particular on the former two, Braess-Sarazin and Vanka, as these are known to lead to effective monolithic multigrid methods, but also expose key kernels that are reused in the implementation of the latter two.  We next provide an overview of all four algorithms, before focusing on aspects of implementation and performance when implementing these approaches on the GPU\@.

    \subsection{Braess-Sarazin relaxation scheme}

        The Braess-Sarazin iteration is based on an approximation of the block factorization of the system matrix in~\eqref{eq:stokesmatrix},
        \begin{equation}
        \begin{bmatrix}
            L & B^T\\
            B & 0
        \end{bmatrix}
        =
                \begin{bmatrix}
            L & 0\\
            B & \hat{S}
        \end{bmatrix}
        \begin{bmatrix}
            I & L^{-1}B^T\\
            0 & I
        \end{bmatrix},
        \end{equation}
        for $\hat{S} = -BL^{-1}B^T$.  The original algorithm~\citep{BraessSarazin} proposed replacing the matrix, $L$, in the above by a scaled version of its diagonal, $tD$, for scalar $t$, and updating the current approximation by an under-relaxed solve of the saddle-point system with this replacement,
        \begin{equation}
            \begin{bmatrix}
                \vec u\\
                p
            \end{bmatrix}^{new}
            =
            \begin{bmatrix}
                \vec u\\
                p
            \end{bmatrix}^{old}
            +
            \omega_{BS}
            \begin{bmatrix}
                tD & B^T\\
                B & 0
            \end{bmatrix}^{-1}
            \begin{bmatrix}
                \vec r_\vec u \\
                r_p
            \end{bmatrix}^{old}
        \end{equation}
        where $\vec r_\vec u$ and $r_p$ are the respective residuals.  The \emph{inexact} variant of Braess-Sarazin~\citep{WZulehner_2000a} computes the (unweighted) updates, $\delta\vec u$ and $\delta p$, as approximate solutions of the block-factorized approximation to the system matrix,
        \begin{equation}
            \begin{bmatrix}
                tD & 0\\
                B & S
            \end{bmatrix}
            \begin{bmatrix}
                I & \frac{1}{t}D^{-1}B^T\\
                0 & I
            \end{bmatrix}
            \begin{bmatrix}
                \delta\vec u\\
                \delta p
            \end{bmatrix}
            =
            \begin{bmatrix}
                \vec {r_u}\\
                r_p
            \end{bmatrix}\label{eq:bsfact}
        \end{equation}
        where $S = -\frac{1}{t}BD^{-1}B^T$ is the Schur complement of the approximated system. \Cref{eq:bsfact} can be rewritten as two equations
        \begin{align}
            S\delta p &= r_p - \frac{1}{t}BD^{-1}\vec{r_u}\label{eq:bseq1}\\
            \delta\vec u &= \frac{1}{t}D^{-1}(\vec{r_u} - B^T\delta p)\label{eq:bseq2}
        \end{align}
        that are solved for $\delta\vec u$ and $\delta p$ using standard weighted Jacobi (or other algorithms) to approximate the inverse of $S$~\citep{WZulehner_2000a,YHe_SMacLachlan_2018b}. The full algorithm is given in~\cref{alg:bs}.
        \begin{algorithm}[ht]
            \begin{algorithmic}[1]
                \State{Approximately solve $S\delta p = r_p - \frac{1}{t}BD^{-1}\vec r_\vec u$ for $\delta p$ by relaxation.}
                \State{Compute $\delta \vec u = \frac{1}{t}D^{-1}(\vec r_\vec u - B^T\delta p)$.}
                \State{Update $p^{new} = p^{old} + \omega_{BS}\delta p$.}
                \State{Update $\vec u^{new} = \vec u^{old} + \omega_{BS}\delta \vec u$.}
            \end{algorithmic}
            \caption{Braess-Sarazin relaxation}\label{alg:bs}
        \end{algorithm}

    \subsection{Vanka relaxation scheme}

\begin{figure}[ht]
            \centering
              \includegraphics{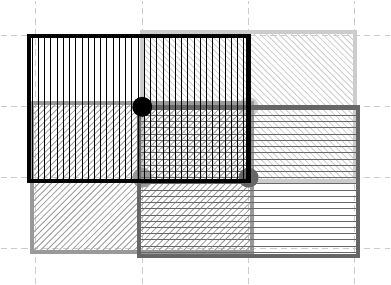}
            \caption{Illustration of overlapping $2\times 2$ Vanka patches}\label{fig:vankapatch}
        \end{figure}

        Vanka relaxation, in contrast, applies a block overlapping Schwarz iteration to the global saddle-point system.  In this approach, we define sets of ``patches'' (or ``subdomains'' in the usual Schwarz notation) corresponding to $2\times 2$ blocks of elements, where we take a single pressure degree of freedom at the central vertex and all velocity degrees of freedom on the associated (neighboring) elements, see~\cref{fig:vankapatch}.  For each patch, we define a restriction operator, $V_i$, that extracts degrees of freedom from the global matrix to those present on local patch $i$, and use this to restrict the system matrix to the $i$th patch, as
        \begin{equation}
            A_i = V_i A V_i^T
        \end{equation}
        The Vanka algorithm is defined by looping over the patches and solving
        \begin{equation}
            A_i\begin{bmatrix}\delta\vec u_i\\\delta p_i\end{bmatrix} = V_i\begin{bmatrix}\vec r_\vec u\\r_p \end{bmatrix}\label{eq:vksystem}
        \end{equation}
        exactly for $\delta \vec u_i$ and $\delta p_i$, which are then used to update the global approximate solution in a weighted additive manner.
        The Vanka algorithm is given in~\cref{alg:vk}. We note that we solve the patch systems exactly by inverting each patch matrix ahead of time, as they do not change between iterations.
        \begin{algorithm}[ht]
            \begin{algorithmic}[1]
                \For{$i = 1, \dots, N$}
                    \State{Solve $A_i\begin{bmatrix}\delta \vec u_i\\\delta p_i\end{bmatrix} = V_i\begin{bmatrix}\vec r_\vec u\\r_p\end{bmatrix}$ for $\begin{bmatrix}\delta \vec u_i\\\delta p_i\end{bmatrix}$.}
                \EndFor
                \State{Update $\begin{bmatrix}\vec u^{new}\\p^{new}\end{bmatrix} = \begin{bmatrix}\vec u^{old}\\p^{old}\end{bmatrix} + \displaystyle\Sigma_{i=1}^NV_i^TW_i\begin{bmatrix}\delta \vec u_i\\\delta p_i\end{bmatrix}$.}
                \Statex{where $W_i$ is the matrix with the weights.}
            \end{algorithmic}
            \caption{Vanka relaxation (additive)}\label{alg:vk}
        \end{algorithm}

    \subsection{Schur-Uzawa relaxation scheme}

        The Schur-Uzawa iteration is derived in a similar way to the Braess-Sarazin iteration. It is again based on an approximation of the factorization of the system matrix in~\eqref{eq:stokesmatrix},
        \begin{equation}
            \begin{bmatrix}
                L & B^T\\
                B & 0
            \end{bmatrix}
            =
                    \begin{bmatrix}
                L & 0\\
                B & \hat{S}
            \end{bmatrix}
            \begin{bmatrix}
                I & L^{-1}B^T\\
                0 & I
            \end{bmatrix},
        \end{equation}
        for $\hat{S} = -BL^{-1}B^T$. As in Braess-Sarazin, we again replace the matrix, $L$, by some approximation that is easy to invert, again denoted $tD$, but also drop the upper-triangular term from this factorization, resulting in an inexact system computing updates, $\delta \vec u$ and $\delta p$, as approximate solutions of the block system,
        \begin{equation}
            \begin{bmatrix}
                tD & 0\\
                B & \hat{S}
            \end{bmatrix}
            \begin{bmatrix}
            \delta \vec u\\
            \delta p
            \end{bmatrix}
            =
            \begin{bmatrix}
                \vec {r_u} \\
                r_p
            \end{bmatrix}.\label{eq:uzblock}
        \end{equation}
        The system in~\eqref{eq:uzblock} can be rewritten as two equations
        \begin{align}
            tD \delta \vec u &= \vec{r_u} \\
            S \delta p &= B \delta\vec u - r_p
        \end{align}
        that are solved for $\delta \vec u$ by directly inverting $\alpha \hat{L}$ and for $\delta p$ by standard weighted Jacobi to approximate the inverse of $S = -\frac{1}{t}BD^{-1}B^T $. The full algorithm is given in~\cref{alg:uz}.
        \begin{algorithm}[ht]
            \begin{algorithmic}[1]
                \State{Compute $\delta \vec u = \frac{1}{t}D^{-1}\vec{r_u} $.}
                \State{Approximately solve $S\delta p = B \delta \vec u - r_p$ for $\delta p$ by relaxation.}
                \State{Update $p^{new} = p^{old} + \delta p$.}
                \State{Update $\vec u^{new} = \vec u^{old} + \delta \vec u$.}
            \end{algorithmic}
            \caption{Schur-Uzawa relaxation}\label{alg:uz}
        \end{algorithm}

    \subsection{Block-Triangular preconditioner}

        The Block-Triangular preconditioner is also based on an approximation of the system matrix in~\eqref{eq:stokesmatrix}, but we now consider an alternate form with unit block diagonal for the lower-triangular factor,
        \begin{equation}
            \begin{bmatrix}
                L & B^T\\
                B & 0
            \end{bmatrix}
            =
            \begin{bmatrix}
                I & 0\\
                BL^{-1} & I
            \end{bmatrix}
            \begin{bmatrix}
                L & B^T\\
                0 & \hat{S}
            \end{bmatrix},
        \end{equation}
        with $\hat{S} = -BL^{-1}B^T$.  While Braess-Sarazin and Schur-Uzawa relaxation use simple approximations to $L$ and $\hat{S}$ to approximate the inverse for relaxation within a multigrid cycle, it is more common to use multigrid on the blocks when using the block preconditioner directly.  Here, as is common~\citep{HCElman_DJSilvester_AWathen_2005a}, we first approximate $\hat{S}$ by a mass matrix, $-M$, on the pressure space, then apply multigrid to this approximation.  With this, we compute updates, $\delta \vec u$ and $\delta p$, as approximate solutions of the upper-triangular approximation to the system matrix,
        \begin{equation}
            \begin{bmatrix}
                L & B^T\\
                0 & -M
            \end{bmatrix}
            \begin{bmatrix}
                \delta \vec u\\
                \delta p
            \end{bmatrix}
            =
            \begin{bmatrix}
            \vec{r_u}\\
            r_p
            \end{bmatrix},\label{eq:schuruzawablock}
        \end{equation}
        using multigrid V-cycles to approximately invert $L$ and $M$.
        The system in~\eqref{eq:schuruzawablock} can be rewritten as two equations
        \begin{align}
            -M \delta p &= r_p\\
            L \delta\vec u &= \vec{r_u} - B^T\delta p,
        \end{align}
leading to the full algorithm given in~\cref{alg:bt}.
        \begin{algorithm}[ht]
            \begin{algorithmic}[1]
                \State{Approximately solve $M\delta p = -r_p$ for $\delta p$ using multigrid on $M$.}
                \State{Approximately solve $L \delta \vec u = \vec{r_u} - B^T\delta p$ for $\delta\vec u$ using multigrid on $L$.}
                \State{Update $p^{new} = p^{old} + \delta p$.}
                \State{Update $\vec u^{new} = \vec u^{old} + \delta \vec u$.}
            \end{algorithmic}
            \caption{Block-Triangular preconditioner}\label{alg:bt}
        \end{algorithm}

\section{Our implementation}

    We have implemented the outer FGMRES iteration and a multigrid V-cycle with the three relaxation schemes, Vanka, Schur-Uzawa, and Braess-Sarazin, as well as the Block-Triangular preconditioner, in C++, with custom data structures that provide structured matrix implementations of the required matrix and vector operations. This is achieved by using operator overloading to allow the optimization of the code for different architectures while preserving a clean implementation of the high-level algorithms. Underlying the custom data structures are standard STL vectors of \lstinline{double} data type. Additionally, support for CUDA and OpenCL requires only a switch of the backend implementation while the high-level algorithm implementation remains largely untouched.

    The resulting implementation yields an efficient solution algorithm for the incompressible Stokes equations in two dimensions on both the CPU and the GPU\@. We limit our attention to optimizing implementations for a single CPU node or single GPU\@ and focus on comparing performance using the different algorithms, taking advantage of the underlying structure. In principle, similar performance is expected for other discretizations of Stokes and other saddle-point problems, in two and three dimensions, but studying performance in these contexts is left for future work.  We also do not consider extending this work to MPI-based parallelism or multi-GPU systems.

\section{Existing work}

    \citet{JohnTobiska} investigate the performance of multigrid paired with Braess-Sarazin for solving the Stokes equations using P1-P0 finite elements. In the case of a W-cycle, the improvement in error reduction is approximately linear with the number of smoothing steps. For a V-cycle, convergence increases in general with increasing level of refinement. The work shows that multigrid paired with Braess-Sarazin is indeed a robust and reliable preconditioner.

    \citet{LarimReusken} compares use of a coupled multigrid method with Vanka and Braess-Sarazin type relaxation schemes, along with preconditioned MINRES and an inexact Uzawa method. The focus is on solving the Stokes equations on the then-current hardware and architectures. The results show that all four methods are robust with respect to variations in parameters. The conclusion is that a multigrid W-cycle paired with diagonal Vanka results in the most efficient solver in terms of CPU time.

    More recently, \citet{JAdler_etal_2015b} compare of a fully-coupled monolithic multigrid paired with Braess-Sarazin or Vanka as relaxation scheme, and a block-factorization preconditioner similar to the one we presented here. On CPU-only systems, multigrid paired with Vanka results in the best scaling and lowest iteration count. Yet, these solvers require significantly more work per iteration than the other preconditioners. As a result, multigrid paired with Braess-Sarazin yields the best time-to-solution on CPUs for the problems studied in that work.

\section{Performance Analysis}

    By far the most costly component of the monolithic multigrid-preconditioned FGMRES solver is the relaxation scheme within the multigrid V-cycle.  Thus, we focus our performance analysis on the implementations of two of the relaxation schemes, Vanka and Braess-Sarazin, alone, noting that multigrid with Schur-Uzawa relaxation reuses only components from that using Braess-Sarazin, while the block-triangular preconditioner also reuses primarily kernels from Braess-Sarazin relaxation as well.

    Studying the performance of these methods requires careful analysis of memory movement and access.  The standard metric for this is \emph{arithmetic intensity}, which quantifies the relationship between floating-point operations and memory reads and writes. Another important quantity is the \emph{FLOP rate}, which describes how many floating-point operations are performed in a given time frame. The runtime of the kernels (and how they relate to one-other) indicates the importance of each kernel when it comes to studying the performance. We first study and optimize the component kernels individually, before comparing performance of our four preconditioners for the FGMRES solver for the Stokes equations.
    For measuring the various metrics to evaluate and compare implementations on the GPU, we use NVIDIA's Nsight Compute\footnote{NVIDIA Nsight Compute: \url{https://developer.nvidia.com/nsight-compute}} and Nsight Systems\footnote{NVIDIA Nsight Systems: \url{https://developer.nvidia.com/nsight-systems}} software.

    \subsection{Test System}

        The system we use for each test is the Delta supercomputer\footnote{Delta supercomputer: \url{https://delta.ncsa.illinois.edu/}} located at the National Center for Supercomputing Applications (NCSA). It is equipped with NVIDIA A100 GPUs that have a measured peak double-precision floating-point performance of 9472.34 GFLOP/s, 80GB on-chip memory, and a measured GPU memory bandwidth of 1264.42 GB/s, measured using the CS roofline toolkit\footnote{CS Roofline Toolkit: \url{https://bitbucket.org/berkeleylab/cs-roofline-toolkit}}. For our final comparison of algorithms, we also run on the CPU nodes of the Delta supercomputer, that carry dual 64-core AMD 7763 processors with a base frequency of 2.45 GHz (max boost frequency of 3.5 GHz) and a per socket memory bandwidth of 204.8 GB/s, although we only consider serial runs on a single core here.

    \subsection{Kernels}

        To start the analysis of the different kernels, we first present an overview of each kernel and get a sense of how much they contribute to the overall runtime.  Although the problem itself and the final parameter choices can have an important impact on the performance of a kernel, we compare kernels for a generic case here to provide a baseline.

\begin{figure}[ht]
            \centering
            \includegraphics[width=\columnwidth]{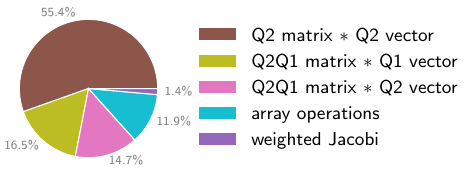}
            \caption{Braess-Sarazin: Kernels and their proportion of runtime}\label{fig:bskernelprop}
        \end{figure}

        \Cref{alg:bs_kernels} shows the Braess-Sarazin algorithm in slightly different form than~\cref{alg:bs}, to focus on the kernels involved for the various steps.  These kernels are color-coded for ease of comparing with the cost breakdown for a single iteration of Braess-Sarazin shown in~\cref{fig:bskernelprop}, indicating how much each kernel contributes to the overall runtime. (Noting that the percentages may not sum to 100\%, due to rounding.)
        As is seen (and expected), most of the runtime is consumed by the matrix-vector operations.  Many of these kernels are reused in the other solvers. We note in particular that the pressure weighted Jacobi kernel, used both for Braess-Sarazin and Schur-Uzawa relaxation (and in the block-triangular preconditioner) contributes very little to the overall runtime (less than $2\%$).
        \begin{algorithm}[ht]
            \begin{algorithmic}[1]
                \State{Compute current residuals, $\vec r_\vec u$ and $r_p$.}
                \Statex{\scriptsize \textbf{\textcolor{tab-brown}{Q2 matrix * Q2 vector}}}
                \Statex{\scriptsize \textbf{\textcolor{tab-pink}{Q2Q1 matrix * Q2 vector}}}
                \Statex{\scriptsize \textbf{\textcolor{tab-olive}{Q2Q1 matrix * Q1 vector}}}
                \Statex{\scriptsize \textbf{\textcolor{tab-cyan}{array operations}}}
                \State{Form right hand side of~\cref{eq:bseq1}.}
                \Statex{\scriptsize \textbf{\textcolor{tab-brown}{Q2 matrix * Q2 vector}}}
                \Statex{\scriptsize \textbf{\textcolor{tab-pink}{Q2Q1 matrix * Q2 vector}}}
                \Statex{\scriptsize \textbf{\textcolor{tab-cyan}{array operations}}}
                \State{Use Jacobi to compute approximation of $\delta p$ in~\cref{eq:bseq1}.}
                \Statex{\scriptsize \textbf{\textcolor{tab-purple}{weighted Jacobi}}}
                \State{Use $\delta p$ to compute $\delta\vec u$ in~\cref{eq:bseq2}.}
                \Statex{\scriptsize \textbf{\textcolor{tab-brown}{Q2 matrix * Q2 vector}}}
                \Statex{\scriptsize \textbf{\textcolor{tab-olive}{Q2Q1 matrix * Q1 vector}}}
                \Statex{\scriptsize \textbf{\textcolor{tab-cyan}{array operations}}}
                \State{Update global solution with $\delta\vec u$ and $\delta p$.}
                \Statex{\scriptsize \textbf{\textcolor{tab-cyan}{array operations}}}
            \end{algorithmic}
            \caption{Braess-Sarazin with kernel breakdowns}\label{alg:bs_kernels}
        \end{algorithm}

\begin{figure*}[ht]
            \centering
            \begin{subfigure}[b]{0.48\linewidth}
                \includegraphics{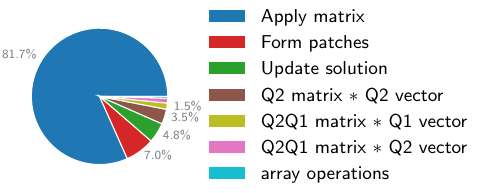}
                \caption{simple Vanka}
            \end{subfigure}
            \begin{subfigure}[b]{0.48\linewidth}
                \includegraphics{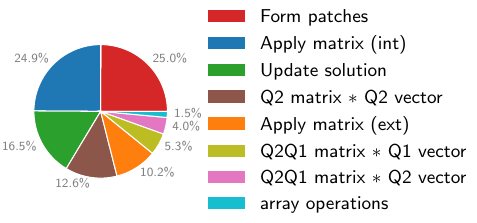}
                \caption{tuned Vanka}
            \end{subfigure}
            \caption{Vanka: Kernels and their proportion of runtime, both (a) simple and (b) tuned Vanka.}\label{fig:vkkernelprop}
        \end{figure*}

\begin{figure}[ht]
            \centering
            \includegraphics{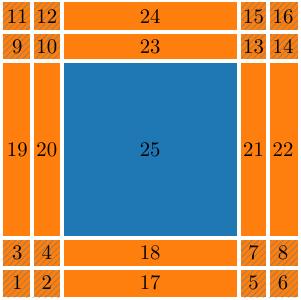}
            \caption{Vanka: groups of shared patch matrices\label{fig:vkmatrices}}
        \end{figure}

        A similar kernel-focused restatement of~\cref{alg:vk} is given in~\cref{alg:vk_kernels}, with the various kernels color-coded to correspond to timing breakdown for a single iteration shown in~\cref{fig:vkkernelprop}.
        \begin{algorithm}[ht]
            \begin{algorithmic}[1]
                \State{Compute current residuals, $\vec r_\vec u$ and $r_p$.}
                \Statex{\scriptsize \textbf{\textcolor{tab-brown}{Q2 matrix * Q2 vector}}}
                \Statex{\scriptsize \textbf{\textcolor{tab-pink}{Q2Q1 matrix * Q2 vector}}}
                \Statex{\scriptsize \textbf{\textcolor{tab-olive}{Q2Q1 matrix * Q1 vector}}}
                \Statex{\scriptsize \textbf{\textcolor{tab-cyan}{array operations}}}
                \State{Form patch right hand sides of~\cref{eq:vksystem}}
                \Statex{\scriptsize \textbf{\textcolor{tab-brown}{Q2 matrix * Q2 vector}}}
                \Statex{\scriptsize \textbf{\textcolor{tab-pink}{Q2Q1 matrix * Q2 vector}}}
                \Statex{\scriptsize \textbf{\textcolor{tab-olive}{Q2Q1 matrix * Q1 vector}}}
                \Statex{\scriptsize \textbf{\textcolor{tab-cyan}{array operations}}}
                \Statex{\scriptsize \textbf{\textcolor{tab-red}{form right hand side}}}
                \State{Apply inverses of patch matrices to patch right hand sides.}
                \Statex{\scriptsize \textbf{\textcolor{tab-blue}{apply matrix inverse (int)}}}
                \Statex{\scriptsize \textbf{\textcolor{tab-orange}{apply matrix inverse (ext)}}}
                \State{Update global solution.}
                \Statex{\scriptsize \textbf{\textcolor{tab-green}{update global solution}}}
            \end{algorithmic}
            \caption{Vanka (tuned) with kernel breakdowns}\label{alg:vk_kernels}
        \end{algorithm}
        At left of~\cref{fig:vkkernelprop}, we show runtime for what we call ``simple Vanka'', where we do not take advantage of the fact that, for many settings, the Vanka submatrices are identical for most patches and can, thus, be stored once and used many times.  This approach results in more than $75\%$ of the runtime being spent applying the patch matrix inverses as each patch needs to load its own Vanka submatrix from global memory. At right, we show results using a ``tuned Vanka'' implementation, where patches that have identical submatrices take advantage of fast shared memory to optimize memory accesses and, in turn, improve performance.  For a uniform grid as considered here, \cref{fig:vkmatrices} sketches the grouping of patches into those that have a submatrix in common.  Here, there are special cases for patches adjacent to the edges or corners of the mesh, including those associated with nodes on the boundary and those distance one from the boundary (where some degrees of freedom in the patch have Dirichlet boundary conditions applied), and a general case for all patches at nodes at least distance two from the boundary. This approach results in only about $40\%$ of the overall runtime being taken up by applying patch matrix inverses. In total, just over three quarters of the runtime in the tuned approach is used for the four unique-to-Vanka operations. The other portion is contributed by the same simple matrix-vector operations as in the analysis of Braess-Sarazin.
        In what follows, we focus on tuned Vanka in our performance analysis and show a comparison of tuned and simple Vanka as part of our final comparison of relaxation schemes within multigrid-preconditioned FGMRES\@.

    \subsection{Vanka patch matrices}

        GPUs are built around multithreaded streaming multiprocessors. Whenever a kernel is launched from the host, all threads are grouped together into smaller thread blocks, which are then enumerated and distributed to available multiprocessors. All threads within a thread block are executed concurrently, and all blocks can be executed concurrently. Threads within a block are able to access local shared memory and can be synchronized. Additionally, all threads are able to access global memory. Choosing the right size of block is essential for good performance, as a too small block size leads to streaming multiprocessors that remain partially idle, whereas a too large block size leads to an imbalanced load over all of the streaming multiprocessors. CUDA is designed with a maximum of 1024 possible threads per thread block.

        Making use of shared and local memory as much as possible within a CUDA thread block allows us to optimize memory accesses further, as data that is used repeatedly can be cached in memory that is faster than global memory. This is of particular importance for the Vanka algorithm, as many of the Vanka patches share the same patch matrices, as discussed above.
        In total, there are $25$ different patch matrices, as depicted in~\cref{fig:vkmatrices}, for a constant-coefficient Stokes problem on a uniform mesh, independent of the number of elements. The orange areas with diagonal lines in~\cref{fig:vkmatrices} denote single patches that have their own unique patch matrix. The normal orange areas are one-dimensional areas along element edges that share the same patch matrix, and the blue area is the two-dimensional interior region of the domain, where all patches share the same submatrix. Within any one of these regions, we can load the patch matrix into shared memory once, to be used by all threads in the block.

    \subsection{Arithmetic Intensity}
        \begin{table*}[t!]
            \small\sf\centering
            \begin{tabular}{m{0.25\textwidth} p{0.2\textwidth}p{0.2\textwidth}p{0.2\textwidth}}
                \toprule
                kernel & reads & writes & flops \\
                \midrule
                array plus/minus array & $2(2n+m)$ & $2n+m$ & $2n+m$ \\ 
                array times scalar & $2n+m$ & $2n+m$ & $2n+m$\\[5pt]
                Q2 matrix $*$ Q2 vector & $n^2+n$ & $n$ & $n^2$\\
                Q2Q1 matrix $*$ Q2 vector & $nm+n$ & $m$ & $nm$\\
                Q2Q1 matrix $*$ Q1 vector & $nm+m$ & $n$ & $nm$\\[5pt]
                Braess-Sarazin: weighted Jacobi & $2m$ & $m$ & $2m$\\
                Vanka: Form Patch RHS & $76 + 124\ell + 51\ell^2$ & $76 + 124\ell + 51\ell^2$ & $0$\\
                Vanka: Apply matrix inverse & $1520 + 3968\ell + 2652\ell^2$ & $76 + 124\ell + 51\ell^2$ & $2888 + 7688\ell + 5202\ell^2$ \\ 
                Vanka: Update global solution & $76 + 124\ell + 51\ell$ & $76 + 124\ell + 51\ell$ & $76 + 124\ell + 51\ell$\\
                \bottomrule
            \end{tabular}
            \caption{Theoretical reads [double], writes [double] and flops of the various operations with $\ell + 2$ as the number of nodal degrees of freedom in one dimension, $n$ as the total number of velocity degrees of freedom, and $m$ as the total number of pressure degrees of freedom.\label{tab:rwf}}
        \end{table*}

        The \emph{arithmetic intensity} of a kernel is defined as the ratio of how many floating point operations (flops) are performed per byte read/written. The algorithms for both Vanka and Braess-Sarazin involve various general vector and matrix-vector operations. In addition, Braess-Sarazin requires a weighted Jacobi application, and Vanka requires operations to extract the current residuals, to apply the patch matrix inverses, and to update the global solution. \Cref{tab:rwf} denotes the counts of all reads, writes, and floating points operations (flops) for the various kernels, obtained by counting the operations in the algorithms.

        Based on the values in~\cref{tab:rwf}, we compute the arithmetic intensity of the various kernels. On the GPU (using CUDA), we use the following formula,
        \begin{equation}
            \text{AI} = \frac{\text{flops}}{32(\text{sectors read $+$ sectors written})}.
        \end{equation}
        where the reads and writes are counted per sector. One sector consists of a total of $32$ bytes and, thus, we multiply by that constant in order to recover the byte count. The performance of each operation is then computed by
        \begin{equation}
            \text{perf} = \frac{\text{flops}}{\max\left(\frac{32(\text{sectors read $+$ written})}{\text{bandwidth}}, \frac{\text{flops}}{\text{peak perf}}\right)}
        \end{equation}
        The theoretical arithmetic intensity and performance computed this way is shown in~\cref{tab:aiperf}.
        \begin{table}
            \centering
            \begin{tabular}{l @{}S@{} @{}S@{}}
                \toprule
                kernel & AI & {performance} \\
                \midrule 
                array plus/minus array    & 0.0417   &  9.821 \\
                array times scalar        & 0.0625   & 14.731   \\[5pt]
                Q2 matrix $*$ Q2 vector   & 0.125    & 29.462 \\
                Q2Q1 matrix $*$ Q2 vector & 0.125    & 29.462 \\
                Q2Q1 matrix $*$ Q1 vector & 0.125    & 29.462 \\[5pt]
                Braess-Sarazin: Jacobi    & 0.0833   & 16.367 \\
                Vanka: Form Patches       & 0.0      &  0.0 \\
                Vanka: Apply matrix       & 0.241    & 56.697 \\
                Vanka: Update solution    & 0.0625   & 14.731   \\
                \bottomrule
            \end{tabular}
            \caption{Theoretical AI [flops/byte] and performance [GFLOP/s], calculated for a $512\times512$ element patch.}\label{tab:aiperf}
        \end{table}
        This analysis, however, has its limitations. In practice, we expect the actual arithmetic intensity and performance to be much better, as values are typically not read from or written to memory one-by-one. Instead, a memory range is typically loaded all at once, allowing us to reuse values. Additional strategies, like using shared memory for Vanka patches with the same patch matrix, further optimize the memory accesses, increasing both the arithmetic intensity and performance. Similarly, varying the size of CUDA blocks also has an effect on these quantities.

    \subsection{Common Kernels}

        For simplicity, we group the kernels into two classes.  First, we examine those kernels that are common to both Vanka and Braess-Sarazin relaxation, involving matrix-vector products and array operations. Following this, we analyze the Vanka-specific kernels.

\begin{figure*}
            \centering
                \begin{tabular}{lll}
                  \mytext{tab-brown}{A}: Q2 matrix $*$ Q2 vector ($n$) &
                  \mytext{tab-olive}{F}: Q2Q1 matrix $*$ Q1 vector ($n$) &
                  \mytext{tab-cyan}{J}: array minus array \\

                  \mytext{tab-brown}{B}: Q2 matrix $*$ Q2 vector ($x$) &
                  \mytext{tab-olive}{G}: Q2Q1 matrix $*$ Q1 vector ($x$) &
                  \mytext{tab-cyan}{K}: array plus array \\

                  \mytext{tab-brown}{C}: Q2 matrix $*$ Q2 vector ($y$) &
                  \mytext{tab-olive}{H}: Q2Q1 matrix $*$ Q1 vector ($y$) &
                  \mytext{tab-cyan}{L}: array plus array (in place)\\

                  \mytext{tab-brown}{D}: Q2 matrix $*$ Q2 vector ($c$) &
                  \mytext{tab-olive}{I}: Q2Q1 matrix $*$ Q1 vector ($c$) &
                  \mytext{tab-cyan}{M}: array times scalar \\

                  \mytext{tab-pink}{E}: Q2Q1 matrix $*$ Q2 vector & & \\
                \end{tabular}\\

                \includegraphics{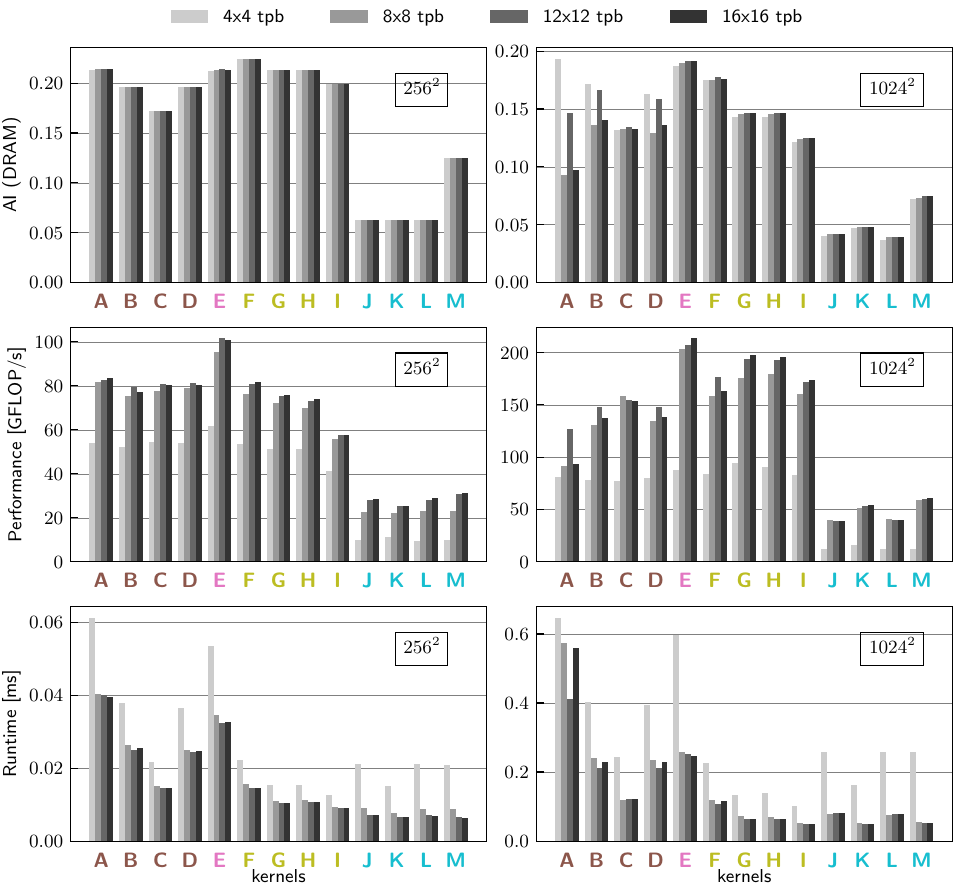}
            \caption{Common kernels: CUDA block size vs. AI, performance, and runtime. $256^2$ elements in left column, $1024^2$ elements in right column.\label{fig:simpletpb}}
        \end{figure*}

        The common kernels are listed at the top of~\cref{fig:simpletpb}, where we break down the matrix-vector products producing Q2 vectors into those that compute values at the nodes, denoted by $n$, the $x$- and $y$-edge midpoints, denoted by $x$ and $y$, respectively, and the cell centers, denoted by $c$.  The color-coding of these kernels matches that in~\cref{fig:bskernelprop}.  For these kernels, we can choose the CUDA block size in an attempt to improve performance. \Cref{fig:simpletpb} shows how the arithmetic intensity, performance, and runtime vary for the various common kernels with varying CUDA block size.

        We first note that the measured arithmetic intensities are indeed better than the theoretical values described in~\cref{tab:aiperf}. Additionally, the measured arithmetic intensity does not vary much (or at all) with varying CUDA block size. This is due to the nature of the underlying memory operations, as the structured matrix data structures used here already optimize the loading and writing of memory.  Due to the global nature of the kernels, they are not able to take advantage of shared memory on the GPU\@. The performance and runtime, however, are impacted by the CUDA block size, with increases in the performance leading to decreases in runtime.  Over all results, we see differences in performance of up to a factor of 5 as we vary the block size.  Choosing the best overall parameter comes down to selecting the best block size for the kernels that contribute the most to each algorithm.

        For Braess-Sarazin, the Q2 matrix by Q2 vector multiplication makes up more than $50\%$ of the overall runtime and, thus, choosing the best parameter for the $4$ kernels within this operation has the largest impact on the overall runtime of the algorithm.  For both problem sizes, the best (or near-best, within $2\%$ of the best) runtime for these $4$ kernels is achieved for a CUDA block size of $12\times12$. Analyzing the other common kernels yields a very similar picture. Thus, all of the common kernels achieve peak (or near-peak) performance for a CUDA block size of $12\times12$, which we choose for the Braess-Sarazin algorithm for which these kernels dominate the cost. We confirmed that these are the best choices by timing a full iteration of the algorithm for both problem sizes. \Cref{tab:bsvkbest} provides a concise overview of the best parameters for both algorithms.

    \subsection{Vanka-Specific Kernels}

\begin{figure*}
            \centering
                  \begin{tabular}{lll}
                    \mytext{tab-orange}{A}: Apply patch matrix (exterior) &
                    \mytext{tab-green}{C}: Update global solution \\
                    \mytext{tab-blue}{B}: Apply patch matrix (interior) &
                    \textbf{\textcolor{tab-red}{D}}: Form patches & \\
                \end{tabular}\\

                \includegraphics{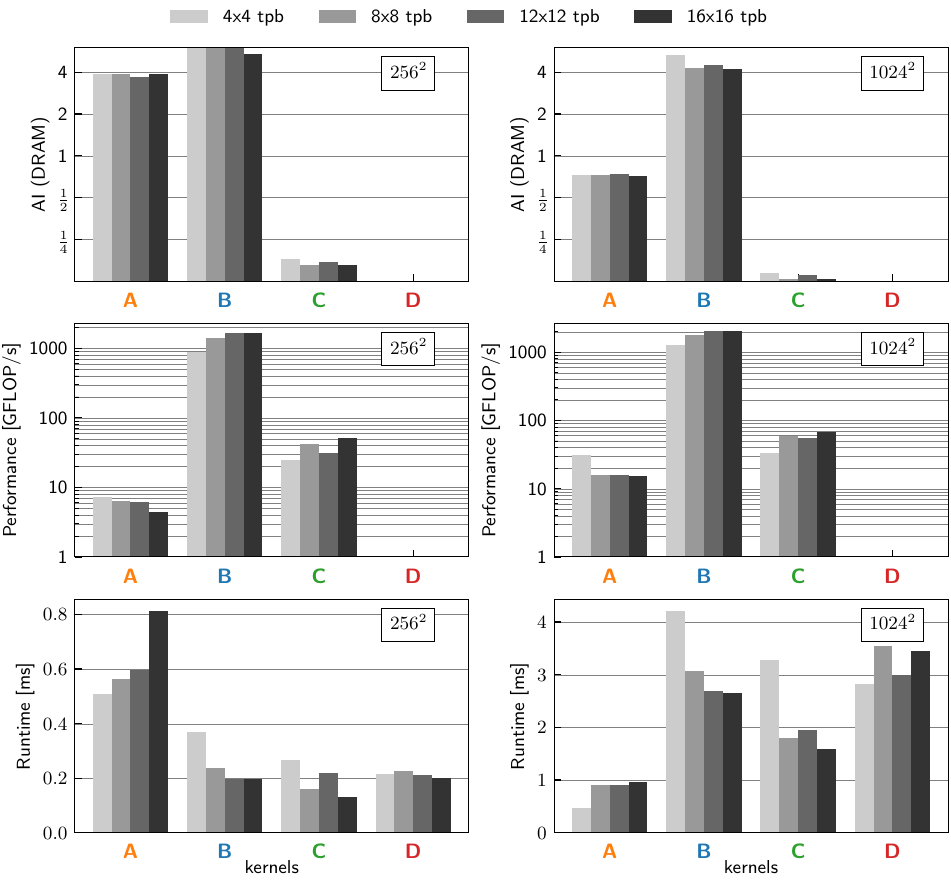}

            \caption{Vanka-specific kernels: thread-block size vs. AI, performance, and runtime. $256^2$ elements in left column, $1024^2$ elements in right column.\label{fig:specifictpb}}
        \end{figure*}

        For the Vanka-specific kernels, we perform a similar analysis as for the common kernels. For all four kernels, we vary the thread block size from $4\times4$ to $16\times16$. \Cref{fig:specifictpb} shows how the arithmetic intensity, performance, and runtime varies with this parameter, again matching the color-coding used in~\cref{fig:vkkernelprop}.
        Here, we notice that the measured arithmetic intensity is higher than the theoretical analysis in~\cref{tab:aiperf}, in particular for the kernels applying the patch inverses. This is expected, as we take advantage of fast shared memory for storing the shared patch matrices, which is not accounted for in that analysis. We note, however, that the arithmetic intensity does not vary much with block size, remaining largely constant. The kernel updating the global solution has a comparatively low arithmetic intensity, as it consists largely of memory movements and only very few floating-point operations. Similarly, the kernel forming the various Vanka patches does not contain any floating-point operations, resulting in zero arithmetic.

        Analyzing the performance of the four kernels shows a rather similar picture, with the thread block size causing only small variations in the performance. Even though the kernel for the exterior patches and the kernel for the interior patches have a very similar arithmetic intensity, they differ widely in terms of performance, by up to $2$ orders of magnitude. This is due to the comparatively high amount of work to be done for the interior patches. Once again, the kernel for forming the Vanka patches has a performance of $0$ GFLOP/s, as it does not contain any floating point operations.

        Both of these metrics, the arithmetic intensity and performance, are useful for evaluating the different kernels, but the effective runtime is the defining criteria for which any set of parameters is, ultimately, the best choice. Even though the kernels applying the patch matrices to the exterior patches (A) has a much lower performance than the kernel applying the patch matrices to the interior patches (B), the runtime of (A) for the smaller problem size is only about a factor of $3$ larger than that for (B). For the larger problem size, the runtime of (A) is much lower than that for (B), by a factor of about $8$. This is due to the overall relatively small amount of computations required for (A), as the exterior regions only grow linearly with the grid size in each dimension, whereas the interior region grows quadratically with (one-dimensional) grid size. Here, we can also see that the proportional runtime for the kernels (B), (C), and (D)  is very much comparable, as already indicated in the kernel runtime breakdown in~\cref{fig:vkkernelprop}.

\begin{figure*}
            \centering
            \includegraphics{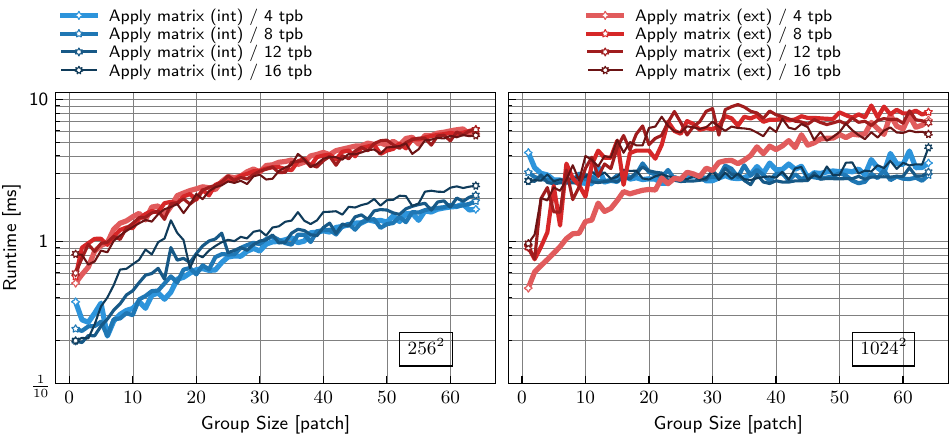}
            \caption{Vanka-specific kernels: Group size vs.\ runtime for $256^2$ elements (left) and $1024^2$ elements (right).\label{fig:specificchunk}}
        \end{figure*}

        Next, we investigate the effect of ``grouping'' computational threads, by passing more than one Vanka patch off to single CUDA thread within any one of the regions where the Vanka patches share the same patch matrix. This reduces the number of overall threads that need to be launched, while potentially further improving the memory accesses required. \Cref{fig:specificchunk} shows the runtime of the two sets of kernels applying the patch matrices for the four different thread block sizes, grouping patches together in groups of $1$ to $64$ patches per thread.
        From~\cref{fig:specificchunk}, we see that increasing the number of patches per thread typically does not lead to a faster runtime; at best, the performance remains relatively constant. Thus, we do not pursue this any further and remain using one thread per patch.

        The four Vanka-specific kernels make up more than $75\%$ of the overall runtime of a Vanka iteration (see~\cref{fig:vkkernelprop}), with each kernel taking up roughly the same proportion of overall runtime. To avoid unnecessary complexity in the code, we choose a single CUDA block size to use for the entire algorithm and all connected kernels. For the smaller problem size, a CUDA block size of $8\times8$ is not the optimal choice for many of the individual kernels, but all four kernels exhibit near-peak performance for this CUDA block size. For the larger problem size, the best choice of CUDA block size is $12\times12$. We have also confirmed that these are the best choices by timing a full iteration of the algorithm for both problem sizes. \Cref{tab:bsvkbest} provides a concise overview of the best parameters for both algorithms.

        \begin{table}[ht]
            \centering
            \begin{tabular}{c c c}
              \toprule
                Algorithm & \# elements & threads/block \\
              \midrule
                \multirow{2}*{Braess-Sarazin} & $256^2$ & $12\times12$ \\
                                              & $1024^2$ & $12\times12$ \\[5pt]
                \multirow{2}*{Vanka} & $256^2$ & $8\times8$\\
                & $1024^2$ & $12\times12$\\
                \bottomrule
            \end{tabular}
            \caption{Best parameter choices for both algorithms.\label{tab:bsvkbest}}
        \end{table}

    \subsection{Roofline model}

\begin{figure*}[ht]
            \centering
            \includegraphics{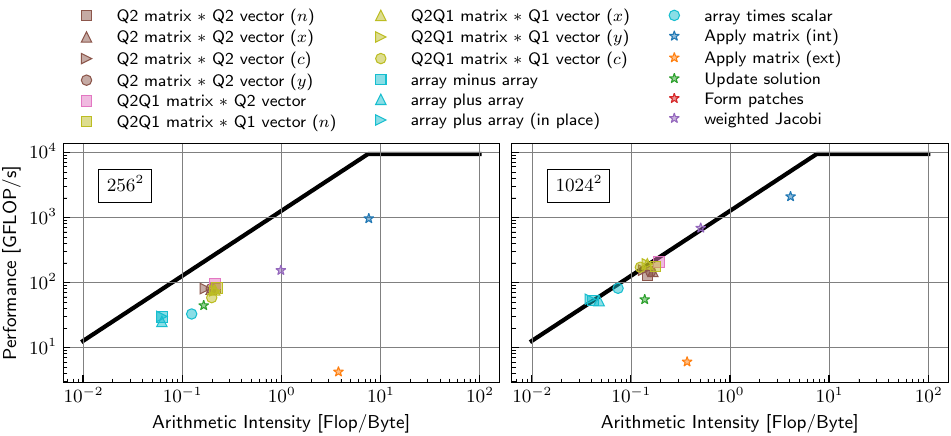}
            \caption{Roofline Model for all kernels}\label{fig:roofline}
        \end{figure*}

        Having analyzed the kernels above and selected the optimal thread block size, we now consider a roofline model to measure for how efficient the kernels are on a given GPU\@. Such models tell us whether an operation is memory or compute bound, and whether all theoretically available computing power is used. \Cref{fig:roofline} shows two roofline models, one for each of the two problem sizes, showing measured performance vs.\ arithmetic intensity for each kernel in a Braess-Sarazin or Vanka relaxation sweep.
        Here, we very clearly see that, for the larger problem size, most kernels lie right on or very close to the performance bound, meaning that they are running as fast as possible given their arithmetic intensity. In order to improve their performance, we would need to find ways to increase their arithmetic intensity. However, given the nature of these kernels and the underlying structured matrix data structures, there is not an obvious avenue to do this.

            Even though most of the kernels (all the matrix-vector and array operations) are clustered together, there are four outliers in this data that we want to highlight:
            \begin{enumerate}
            \item The first outlier is the kernel corresponding to forming the Vanka patch submatrices, which is not visible in the plot, as it consists entirely of memory movements and no floating-point operations. Its arithmetic intensity and floating-point performance are, thus, $0$.
                \item The second outlier is the kernel applying the patch matrix inverse to the exterior patches of the domain. This has a low peak performance of only $6$ GFLOP/s, as it consists of mostly small operations (16 unique patch matrices, with 8 one-dimensional regions sharing a patch matrix).  It also acts on little enough data that, even for the large problem size, we do not achieve the performance expected from the roofline model.
                \item The third outlier is a kernel that we mostly ignored in our analysis, the weighted Jacobi kernel. This kernel achieves a higher performance than all but one other kernel, with a peak performance of $753$ GFLOP/s. However, it contributes less than $2\%$ to the overall runtime of Braess-Sarazin relaxation, with similar percentages of runtime expected for the other algorithms that use it. Thus, even though its performance is rather high, it has barely any measurable effect on the algorithm runtime.
                \item The final outlier is the kernel applying the patch matrix inverse to the interior patches of the domain. Its peak performance is roughly $2145$ GFLOP/s, almost three times as high as the next highest kernel. With this high performance, it still makes up about $20\%$ of the overall runtime of Vanka. Thus, achieving this performance on this single kernel results in the overall Vanka algorithm achieving much better performance.
            \end{enumerate}

            Overall, we note that most of the kernels achieve their maximum possible performance, as they lie right on the performance limit in the roofline model for the larger problem size. Due to the nature of their operations, increasing their arithmetic intensity is not possible and, thus, the performance of these algorithms cannot reasonably be expected to be increased. One avenue to consider to improve the performance of the kernels that require a disproportionately large volume of memory movement would be to try to ``trade'' some memory movement for increasing numbers of floating-point operations; this will be a subject for future research.

    \subsection{Overall solver performance}

        Finally, having analyzed and optimized the performance of the Vanka and Braess-Sarazin relaxation schemes, we now look to see how they compare in practice, when used as relaxation schemes inside of a multigrid V-cycle that is used as preconditioner for FGMRES applied to the Stokes equations. We will also compare their performances to the performance of FGMRES preconditioned with a multigrid V-cycle with Schur-Uzawa and preconditioned with a Block-Triangular preconditioner with multigrid V-cycles used to approximate the block inverses. The additional parameters needed for Schur-Uzawa and the Block-Triangular preconditioner have been determined through further experiment. The optimal value of $t$ in the Schur-complement scheme is $1$ with an optimal Jacobi weight of $\omega=0.4$. For the Block-Triangular preconditioner, we determined that a total of $3$ V-cycles are necessary for both the pressure update solve and velocity update solve, and the two respective weights for the weighted Jacobi relaxation are $\omega=0.6$ for the pressure update solve, and $\omega = 1.0$ for the velocity update solve.

        We use our own implementation of FGMRES, making use of our structured data structures, and use a multigrid \lstinline{V(1,1)} cycle as preconditioner, and a \lstinline{V(3,3)} cycle as part of the Block-Triangular preconditioner. At each level of the multigrid algorithm, we use a sweep of either Braess-Sarazin, Vanka, or Schur-Uzawa relaxation. With the Block-Triangular preconditioner we use three sweeps of weighted Jacobi. At the coarsest level, we use either an exact solve on the CPU or three sweeps of the relaxation scheme on the GPU\@.

\begin{figure}[ht]
            \centering
            \includegraphics{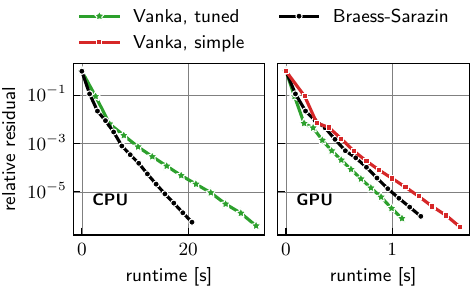}
            \caption{Comparing Vanka to Braess-Sarazin for a problem of size $1024^2$ elements ($768^2$ elements for simple Vanka)}\label{fig:compvkbs}
        \end{figure}

        The first comparison we consider is a comparison of Braess-Sarazin to both our tuned Vanka implementation described in this paper and a simple Vanka implementation, shown in \cref{fig:compvkbs}.
        All of the results are for problems of size $1024^2$, with the exception of the simple Vanka runs. Due to its higher memory requirements, the largest problem size that successfully ran was a problem of size $768^2$ elements. However, even though simple Vanka has just over half as many elements as the other approaches, it is still not able match their performance. On the CPU, we see that multigrid with Braess-Sarazin relaxation strongly outperforms the use of Vanka relaxation. Even though multigrid with Vanka relaxation typically requires one fewer iteration to reach convergence, the work required per iteration is significantly larger than for Braess-Sarazin, resulting in multigrid with Vanka taking about twice as long. On the GPU, however, we are able to take advantage of the throughput of Vanka, resulting in a runtime that is more than $23$ times smaller than on the CPU, whereas the runtime for Braess-Sarazin is only reduced by a factor of about $11$. Overall, on the GPU, tuned Vanka outperforms Braess-Sarazin by about $10\%$.

\begin{figure}[ht]
            \centering
            \includegraphics{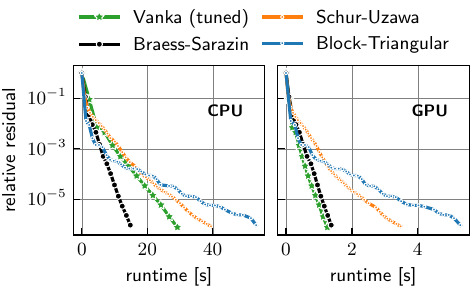}
            \caption{Comparing all four preconditioning strategies for a problem of size $1024^2$ elements}\label{fig:compall}
        \end{figure}

        Next, we compare Braess-Sarazin and tuned Vanka to the other two preconditioning strategies, monolithic multigrid with Schur-Uzawa and the block-triangular preconditioner, shown in \cref{fig:compall}.
        We can see that both the multigrid preconditioner with Schur-Uzawa relaxation and the Block-Triangular preconditioner are not able to match the performance of both Braess-Sarazin and Vanka. Initially they perform very well, in particular the Block-Triangular preconditioner, but they soon slow down requiring up to more than 3 times as long as Braess-Sarazin and Vanka (on the GPU).

\begin{figure}[ht]
            \centering
            \includegraphics{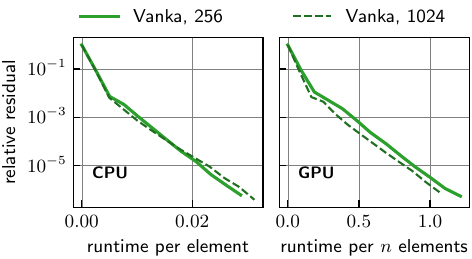}
            \caption{Showing the work per element (in ms) of tuned Vanka for $256^2$ and $1024^2$ elements.}\label{fig:compscaling}
        \end{figure}

        The third metric to consider is the performance of our tuned Vanka implementation for two problem sizes when normalized per element on the CPU and per row of elements on the GPU\@; this is shown in \cref{fig:compscaling}.
        We observe that the time for the tuned Vanka implementation (per element) remains the same no matter the problem size on the CPU, requiring about $0.03$ms per element. Thus, there is no additional overhead introduced by the size of the problem. On the GPU, we are able expose the fine-grained parallelism in the tuned Vanka implementation, resulting in a constant work per row of elements at just over $1$ ms. In fact, we are able to perform about $10\%$ faster per row of elements for the larger problem size.

        These results show that a careful implementation of Vanka on the GPU not only results in the fastest time to convergence, but it also does so without requiring additional parameters to be set. In addition, the parallelism of Vanka makes it a clear favorite in distributed memory settings.

\section{Conclusions and Future Work}

    Several preconditioners for FGMRES are well-known to yield scalable solution algorithms for saddle-point problems, such as the Stokes equations, including both monolithic multigrid and multigrid-based block-factorization preconditioners.  Here, we consider their implementation, performance, and optimization, on modern CPU and GPU architectures.  Different metrics were presented and analyzed, including arithmetic intensity, performance, and runtime, for the various kernels making up these algorithms. Given a highly structured setup, we show that multigrid with Vanka relaxation can be very performant on the GPU\@, leading to faster convergence than when using Braess-Sarazin, both in terms of iterations (saving just $1$ iteration) and runtime (up to $10\%$ faster). This shows that using Vanka relaxation is both mathematically and computationally competitive, although a careful design of the algorithm is warranted. This also highlights the benefit of using GPUs for such algorithms, as multigrid with Vanka on the GPU is up to $23$ times faster than on the CPU, while multigrid with Braess-Sarazin is up to $11$ times faster.

    We also presented two other preconditioning strategies, multigrid preconditioner with Schur-Uzawa relaxation, and Block-Triangular preconditioner with multigrid and weighted Jacobi within. Both of these have been shown to not be able to compete with multigrid with Braess-Sarazin or Vanka, in particular on the GPU\@. In addition, they introduce additional parameters that need to be carefully chosen.

    Future work includes extending this work to cases where the tuned Vanka approach is not applicable, such as for linearizations of the Navier-Stokes equations.  It is also not clear how well these results generalize to other discretizations of saddle-point systems (including both higher-order discretizations using generalized Taylor-Hood elements and other discretizations, such as using discontinuous Galerkin methods).  Extensions to three-dimensional incompressible flow problems and other saddle-point systems are also important future work.  One such system of interest, for example, that combines some of these difficulties is the Reynolds-Averaged Navier-Stokes (RANS) equations, in the context of wind-turbine simulations.

\section{Acknowledgments}

    This work used the Delta system at the National Center for Supercomputing Applications through allocation CIS230037 from the Advanced Cyberinfrastructure Coordination Ecosystem: Services \& Support (ACCESS) program, which is supported by National Science Foundation grants \#2138259, \#2138286, \#2138307, \#2137603, and \#2138296.

\section{Resources}

    The implementation of FGMRES and multigrid with Braess-Sarazin, Vanka, and Schur-Uzawa relaxation, and the Block-Triangular preconditioner, as well as problem setup and all data structures is available at: \url{https://gitlab.com/luspi/gmres}.

\bibliographystyle{SageH}

\begin{thebibliography}{34}
\providecommand{\natexlab}[1]{#1}
\providecommand{\url}[1]{\texttt{#1}}
\providecommand{\urlprefix}{URL }
\expandafter\ifx\csname urlstyle\endcsname\relax
  \providecommand{\doi}[1]{DOI:\discretionary{}{}{}#1}\else
  \providecommand{\doi}{DOI:\discretionary{}{}{}\begingroup
  \urlstyle{rm}\Url}\fi

\bibitem[{Adler et~al.(2023)Adler, He, Hu, MacLachlan and
  Ohm}]{JAdler_etal_2023a}
Adler J, He Y, Hu X, MacLachlan S and Ohm P (2023) Monolithic multigrid for a
  reduced- quadrature discretization of poroelasticity.
\newblock \emph{SIAM Journal on Scientific Computing} 45(3): S54--S81.

\bibitem[{Adler et~al.(2016)Adler, Benson, Cyr, MacLachlan and
  Tuminaro}]{BensonAdlerCyrMacLachlanTuminaro}
Adler JH, Benson TR, Cyr EC, MacLachlan SP and Tuminaro RS (2016) Monolithic
  multigrid methods for two-dimensional resistive magnetohydrodynamics.
\newblock \emph{SIAM Journal on Scientific Computing} 38(1): B1--B24.
\newblock \doi{10.1137/151006135}.
\newblock \urlprefix\url{http://dx.doi.org/10.1137/151006135}.

\bibitem[{Adler et~al.(2017)Adler, Benson and MacLachlan}]{JAdler_etal_2015b}
Adler JH, Benson TR and MacLachlan SP (2017) Preconditioning a mass-conserving
  discontinuous {G}alerkin discretization of the {S}tokes equations.
\newblock \emph{Numerical Linear Algebra with Applications} 24(3): e2047.
\newblock \doi{10.1002/nla.2047}.

\bibitem[{Ayuso~de Dios et~al.(2014)Ayuso~de Dios, Brezzi, Marini, Xu and
  Zikatanov}]{BAyuso_etal_2014a}
Ayuso~de Dios B, Brezzi F, Marini LD, Xu J and Zikatanov L (2014) A simple
  preconditioner for a discontinuous {G}alerkin method for the {S}tokes
  problem.
\newblock \emph{Journal of Scientific Computing} 58(3): 517--547.
\newblock \doi{10.1007/s10915-013-9758-0}.
\newblock \urlprefix\url{http://dx.doi.org/10.1007/s10915-013-9758-0}.

\bibitem[{Benzi et~al.(2005)Benzi, Golub and
  Liesen}]{MBenzi_GHGolub_JLiesen_2005a}
Benzi M, Golub G and Liesen J (2005) Numerical solution of saddle point
  problems.
\newblock \emph{Acta Numerica} 14: 1--137.

\bibitem[{Bienz et~al.(2016)Bienz, Falgout, Gropp, Olson and
  Schroder}]{2016_BiFaGrOlSc}
Bienz A, Falgout RD, Gropp W, Olson LN and Schroder JB (2016) Reducing parallel
  communication in algebraic multigrid through sparsification.
\newblock \emph{SIAM Journal on Scientific Computing} 38(5): S332--S357.

\bibitem[{Bienz et~al.(2020)Bienz, Gropp and Olson}]{2020_BiGrOl_reducing}
Bienz A, Gropp WD and Olson LN (2020) Reducing communication in algebraic
  multigrid with multi-step node aware communication.
\newblock \emph{The International Journal of High Performance Computing
  Applications} 34(5): 547--561.

\bibitem[{Braess and Sarazin(1997)}]{BraessSarazin}
Braess D and Sarazin R (1997) An efficient smoother for the {S}tokes problem.
\newblock \emph{Applied Numerical Mathematics} 23(1): 3--19.
\newblock \doi{10.1016/S0168-9274(96)00059-1}.
\newblock \urlprefix\url{http://dx.doi.org/10.1016/S0168-9274(96)00059-1}.
\newblock Multilevel methods (Oberwolfach, 1995).

\bibitem[{Brandt and Dinar(1979)}]{MultigridStokes}
Brandt A and Dinar N (1979) Multigrid solutions to elliptic flow problems.
\newblock In: PARTER SV (ed.) \emph{Numerical Methods for Partial Differential
  Equations}. Academic Press.
\newblock ISBN 978-0-12-546050-7, pp. 53--147.
\newblock \doi{https://doi.org/10.1016/B978-0-12-546050-7.50008-3}.
\newblock
  \urlprefix\url{https://www.sciencedirect.com/science/article/pii/B9780125460507500083}.

\bibitem[{Briggs et~al.(2000)Briggs, Henson and
  McCormick}]{WLBriggs_VEHenson_SFMcCormick_2000a}
Briggs WL, Henson VE and McCormick SF (2000) \emph{A Multigrid Tutorial}.
\newblock Philadelphia: SIAM Books.
\newblock Second edition.

\bibitem[{Dendy(1982)}]{JEDendy_1982a}
Dendy JE (1982) Black box multigrid.
\newblock \emph{J. Comput. Phys.} 48: 366--386.

\bibitem[{Dou and Liang(2023)}]{SolverCurrent1}
Dou Y and Liang ZZ (2023) A class of block alternating splitting implicit
  iteration methods for double saddle point linear systems.
\newblock \emph{Numerical Linear Algebra with Applications} 30(1): e2455.
\newblock \doi{https://doi.org/10.1002/nla.2455}.
\newblock
  \urlprefix\url{https://onlinelibrary.wiley.com/doi/abs/10.1002/nla.2455}.

\bibitem[{Elman et~al.(2005)Elman, Silvester and
  Wathen}]{HCElman_DJSilvester_AWathen_2005a}
Elman H, Silvester D and Wathen A (2005) \emph{Finite elements and fast
  iterative solvers: with applications in incompressible fluid dynamics}.
\newblock Numerical Mathematics and Scientific Computation. New York: Oxford
  University Press.
\newblock ISBN 978-0-19-852868-5; 0-19-852868-X.

\bibitem[{Ershkov et~al.(2021)Ershkov, Prosviryakov, Burmasheva and
  Christianto}]{SolverCurrent3}
Ershkov SV, Prosviryakov EY, Burmasheva NV and Christianto V (2021) Towards
  understanding the algorithms for solving the {N}avier–{S}tokes equations.
\newblock \emph{Fluid Dynamics Research} 53(4): 044501.
\newblock \doi{10.1088/1873-7005/ac10f0}.
\newblock \urlprefix\url{https://dx.doi.org/10.1088/1873-7005/ac10f0}.

\bibitem[{Farrell et~al.(2021)Farrell, He and MacLachlan}]{FourierVankaStokes}
Farrell PE, He Y and MacLachlan SP (2021) A local {F}ourier analysis of
  additive {V}anka relaxation for the {S}tokes equations.
\newblock \emph{Numerical Linear Algebra with Applications} 28(3): e2306.
\newblock \doi{https://doi.org/10.1002/nla.2306}.
\newblock
  \urlprefix\url{https://onlinelibrary.wiley.com/doi/abs/10.1002/nla.2306}.

\bibitem[{Greif and He(2021)}]{VankaSmoothing}
Greif C and He Y (2021) A closed-form multigrid smoothing factor for an
  additive {V}anka-type smoother applied to the {P}oisson equation.

\bibitem[{He and MacLachlan(2019)}]{YHe_SMacLachlan_2018b}
He Y and MacLachlan S (2019) Local {F}ourier analysis for mixed finite-element
  methods for the {S}tokes equations.
\newblock \emph{Journal of Computational and Applied Mathematics} 357:
  161--183.

\bibitem[{He and MacLachlan(2020)}]{YHe_SMacLachlan_2018a}
He Y and MacLachlan S (2020) Two-level {F}ourier analysis of multigrid for
  higher-order finite-element discretizations of the {L}aplacian.
\newblock \emph{Numer. Linear Alg. Appl.} 27(3): e2285.

\bibitem[{John and Tobiska(2000{\natexlab{a}})}]{JohnTobiska}
John V and Tobiska L (2000{\natexlab{a}}) A coupled multigrid method for
  nonconforming finite element discretizations of the 2{D}-{S}tokes equation.
\newblock \emph{Computing} 64(4): 307--321.
\newblock \doi{10.1007/s006070070027}.
\newblock \urlprefix\url{http://dx.doi.org/10.1007/s006070070027}.
\newblock International GAMM-Workshop on Multigrid Methods (Bonn, 1998).

\bibitem[{John and Tobiska(2000{\natexlab{b}})}]{VKBS1}
John V and Tobiska L (2000{\natexlab{b}}) Numerical performance of smoothers in
  coupled multigrid methods for the parallel solution of the incompressible
  {N}avier-{S}tokes equations.
\newblock \emph{International Journal For Numerical Methods In Fluids} 33(4):
  453--473.

\bibitem[{Larin and Reusken(2008{\natexlab{a}})}]{CompVankaBSStokes}
Larin M and Reusken A (2008{\natexlab{a}}) A comparative study of efficient
  iterative solvers for generalized {S}tokes equations.
\newblock \emph{Numerical Linear Algebra with Applications} 15(1): 13--34.
\newblock \doi{10.1002/nla.561}.
\newblock \urlprefix\url{http://dx.doi.org/10.1002/nla.561}.

\bibitem[{Larin and Reusken(2008{\natexlab{b}})}]{LarimReusken}
Larin M and Reusken A (2008{\natexlab{b}}) A comparative study of efficient
  iterative solvers for generalized {S}tokes equations.
\newblock \emph{Numerical Linear Algebra with Applications} 15(1): 13--34.
\newblock \doi{10.1002/nla.561}.
\newblock \urlprefix\url{http://dx.doi.org/10.1002/nla.561}.

\bibitem[{Maitre et~al.(1984)Maitre, Musy and
  Nignon}]{JFMaitre_FMusy_PNignon_1984a}
Maitre JF, Musy F and Nignon P (1984) A fast solver for the {S}tokes equations
  using multigrid with a {UZAWA} smoother.
\newblock In: Braess D, Hackbusch W and Trottenberg U (eds.) \emph{Advances in
  Multi--Grid Methods}, \emph{Notes on Numerical Fluid Mechanics}, volume~11.
  Braunschweig: Vieweg, pp. 77--83.

\bibitem[{Munch and Kronbichler(2023)}]{MunchKronbichlerAS}
Munch P and Kronbichler M (2023) Cache-optimized and low-overhead
  implementations of additive schwarz methods for high-order fem multigrid
  computations.
\newblock \emph{The International Journal of High Performance Computing
  Applications} : 10943420231217221\doi{10.1177/10943420231217221}.
\newblock \urlprefix\url{https://doi.org/10.1177/10943420231217221}.

\bibitem[{Nataf and Tournier(2022)}]{SolverCurrent2}
Nataf F and Tournier PH (2022) Recent advances in domain decomposition methods
  for large-scale saddle point problems.
\newblock \emph{Comptes Rendus. M\'ecanique} \doi{10.5802/crmeca.130}.
\newblock Online first.

\bibitem[{Notay(2019)}]{SolverBlock}
Notay Y (2019) Convergence of some iterative methods for symmetric saddle point
  linear systems.
\newblock \emph{SIAM Journal on Matrix Analysis and Applications} 40(1):
  122--146.
\newblock \doi{10.1137/18M1208836}.
\newblock \urlprefix\url{https://doi.org/10.1137/18M1208836}.

\bibitem[{Paisley and Bhatti(1998)}]{VKBS2}
Paisley M and Bhatti N (1998) Comparison of multigrid methods for neutral and
  stably stratified flows over two-dimensional obstacles.
\newblock \emph{Journal of Computational Physics} 142(2): 581–610.
\newblock \doi{10.1006/jcph.1998.5915}.
\newblock \urlprefix\url{https://doi.org/10.1006/jcph.1998.5915}.

\bibitem[{Reisner et~al.(2020)Reisner, Berndt, Moulton and
  Olson}]{ScalableLineRelaxation}
Reisner A, Berndt M, Moulton JD and Olson LN (2020) Scalable line and plane
  relaxation in a parallel structured multigrid solver.
\newblock \emph{Parallel Computing} 100: 102705.
\newblock \doi{https://doi.org/10.1016/j.parco.2020.102705}.
\newblock
  \urlprefix\url{https://www.sciencedirect.com/science/article/pii/S0167819120300922}.

\bibitem[{Reisner et~al.(2018)Reisner, Olson and Moulton}]{ScalingStructuredMG}
Reisner A, Olson LN and Moulton JD (2018) Scaling structured multigrid to
  500{K}+ cores through coarse-grid redistribution.
\newblock \emph{SIAM Journal on Scientific Computing} 40(4): C581--C604.
\newblock \doi{10.1137/17M1146440}.
\newblock \urlprefix\url{https://doi.org/10.1137/17M1146440}.

\bibitem[{Trottenberg et~al.(2001)Trottenberg, Oosterlee and
  Sch{\"u}ller}]{UTrottenberg_etal_2001a}
Trottenberg U, Oosterlee CW and Sch{\"u}ller A (2001) \emph{Multigrid}.
\newblock London: Academic Press.

\bibitem[{ur~Rehman et~al.(2011)ur~Rehman, Geenen, Vuik, Segal and
  MacLachlan}]{StokesIterative}
ur~Rehman M, Geenen T, Vuik C, Segal G and MacLachlan SP (2011) On iterative
  methods for the incompressible {S}tokes problem.
\newblock \emph{International Journal for Numerical Methods in Fluids} 65(10):
  1180--1200.
\newblock \doi{https://doi.org/10.1002/fld.2235}.
\newblock
  \urlprefix\url{https://onlinelibrary.wiley.com/doi/abs/10.1002/fld.2235}.

\bibitem[{Vanka(1986)}]{Vanka2}
Vanka S (1986) Block-implicit multigrid solution of {N}avier-{S}tokes equations
  in primitive variables.
\newblock \emph{Journal of Computational Physics} 65(1): 138--158.
\newblock \doi{https://doi.org/10.1016/0021-9991(86)90008-2}.
\newblock
  \urlprefix\url{https://www.sciencedirect.com/science/article/pii/0021999186900082}.

\bibitem[{Voronin et~al.(2022)Voronin, He, MacLachlan, Olson and
  Tuminaro}]{LowOrderPreconditioning}
Voronin A, He Y, MacLachlan S, Olson LN and Tuminaro R (2022) Low-order
  preconditioning of the {S}tokes equations.
\newblock \emph{Numerical Linear Algebra with Applications} 29(3): e2426.
\newblock \doi{https://doi.org/10.1002/nla.2426}.
\newblock
  \urlprefix\url{https://onlinelibrary.wiley.com/doi/abs/10.1002/nla.2426}.

\bibitem[{Zulehner(2000)}]{WZulehner_2000a}
Zulehner W (2000) A class of smoothers for saddle point problems.
\newblock \emph{Computing} 65(3): 227--246.

\end{thebibliography}

\end{document}